\documentclass[letterpaper,landscape]{article}
\usepackage{mathptmx,graphicx}
\twocolumn
\def\tab{\hskip .166667in}
\def\tabb{\hskip .666667in}
\def\<{\langle}
\def\>{\rangle}

\begin{document}
\phantom{}

\eject{\footnotesize Discrete Comput. Geom. {\bf 40}(2), 214-240 (2008)
\par DOI 10.1007/s00454-008-9101-y
\par Published online: 22 August 2008

}\vskip .5in{\large\bf A Dense Packing of Regular Tetrahedra}

\bigskip
\bigskip{\bf Elizabeth R. Chen}

\vskip 1in{\footnotesize  Received: 1 January 2007 / Revised: 9 July 2008 / Accepted: 9 July 2008
\par\copyright\, Springer Science+Business Media, LLC

\bigskip}
\bigskip{\bf Abstract}\, We construct a dense packing of regular tetrahedra, with packing density $D > .7786157$.


\bigskip
\bigskip{\bf 1\, Introduction}

\bigskip ``How can one arrange most densely in space an infinite number of equal solids of given form, e.g., spheres with given radii or regular tetrahedra with given edges (or in prescribed position), that is, how can one so fit them together that the ratio of the filled to the unfilled space may be as great as possible?'' (excerpt from David Hilbert's 18th problem [8]).

\par\tab This is a very old problem. Aristotle [1] believed that you could tile space with regular tetrahedra. Everyone believed him for the next 1800 years, until Johannes M\"uller (aka Regimontanus) proved everyone wrong (histories by Dirk Struik [13] and Marjorie Senechal [12]).

\par\tab This is a very understandable mistake, because tetrahedra almost tile space locally, but not quite. The cluster ${\bf E}_5^{}$ (consisting of 5 tetrahedra joined symmetrically about an edge) has total solid angle $5\cdot 2\cos^{-1}{1\over 3}\approx 12.309594173408$ about the edge, and local density $D\approx .979566380077$. The cluster ${\bf V}_{20}^{}$ (consisting of 20 tetrahedra joined symmetrically about a vertex) has total solid angle $20\cdot(-\pi+3\cos^{-1}{1\over 3})\approx 11.025711968651$ about the vertex, and local density $D\approx .877398280459$.

\bigskip
\bigskip{\footnotesize This research was partially supported by the NSF-RTG grant DMS-0502170.

\bigskip E.R. Chen (corresponding author)
\par Department of Mathematics, University of Michigan, Ann Arbor, MI 48109, USA
\par e-mail: bethchen@umich.edu

}\eject\tab For the sphere, the analogous problem was very challenging. A long time ago, Johannes Kepler conjectured that the densest packing is the hexagonal close-packing (HCP), with density $D = \pi/\sqrt{18}\approx .740480489693$. Carl Friedrich Gauss proved that the HCP is the densest lattice packing. Only recently, Thomas Hales and Samuel Ferguson [6, 7] proved that the HCP is the densest packing in general.

\par\tab Helmut Gr\"omer [5] constructed a lattice packing of the single tetrahedron ${\bf B}_1^{}$ with density $D = {18\over 49}\approx .367346938775$. Douglas Hoylman [9] proved that Gr\"omer's packing was the densest lattice packing.

\par\tab Andrew Hurley [10] constructed the tetrahelix. I calculated its packing density $D = \sqrt{50000\over 177417}\approx .531273435694$.

\par\tab Hermann Minkowski [11] proved that the densest lattice packing of any convex body must satisfy certain constraints. Ulrich Betke \& Martin Henk [2, 3] developed an efficient computer algorithm to compute the densest lattice packing of any convex body, and they applied it to the Archimedean solids. I used Betke \& Henk's program to calculate the packing density of the convex hull of ${\bf V}_{20}^{}$, $D\approx .716796401602$.

\par\tab John Conway \& Salvatore Torquato [4] used Betke \& Henk's program to calculate the packing density of ${\bf V}_{20}^{}$ inscribed in a regular icosahedron, $D > .7165598$. They slightly improved the packing density by wiggling the tetrahedra, $D > .717455$. They raised the question whether the packing density of regular tetrahedra is, not only less than that of the sphere, but perhaps even the lowest of any convex body.

\par\tab In this paper, we will construct a dense packing of regular tetrahedra
with packing density $D\approx .778615700855$, which beats the densest sphere packing!

\bigskip
\bigskip{\bf 2\, Description by Construction}

\bigskip Feel free to look at the pictures in Sec. 3, as you read the descriptions in Sec. 2.

\par\tab We construct a 2-parameter family of clusters of regular tetrahedra. Each cluster (9 tetrahedra) is the union of 1 `central' tetrahedron, 4 `upper' tetrahedra attached to an edge of the central tetrahedron, and 4 `lower' tetrahedra attached to the opposite edge of the central tetrahedron. The 4 upper tetrahedra can rotate about its attached edge, by an angle parametrized by $u$, and the 4 lower tetrahedra can rotate about its attached edge, by an angle parametrized by $v$. Clusters have 2 orientations: `positive' and `negative' (scalar mult. by $-1$).

\par\tab For each parameter value $\<u,v\>$, we construct an optimized packing of clusters, which is crystallographic. The clusters pack in layers with alternating orientations. The full symmetry group of the packing acts transitively on all clusters. The direct symmetry group, which is also the translation symmetry group, acts transitively on a coset. There are 2 cosets, which correspond to the 2 orientations = 2 translation classes. So the fundamental domain (under the lattice of translations) contains 2 clusters = 18 tetrahedra.

\par\tab Since every cluster is equivalent to all other clusters, we restrict our attention to the cluster at the origin and its immediate neighbors. Also, since the cluster (9 tetrahedra) is extremely non-convex, we think of it as the union of (the convex hull of) the `upper half-cluster' (5 tetrahedra) and (the convex hull of) the `lower half-cluster' (5 tetrahedra). In order to check that clusters don't overlap in the packing, it's convenient to check that half-clusters don't overlap.

\eject\tab Each cluster has transverse edge-to-edge intersections with neighboring clusters in the same layer, and partial face-to-face intersections with neighboring clusters in adjacent layers. For intersecting edges, the separating plane contains the union of the edges, and for non-intersecting edges, the separating plane is between and parallel to both edges.
For intersecting faces, the separating plane contains the union of the faces, and for non-intersecting faces, the separating plane is between and parallel to both faces.

\par\tab The intersections determine equations in terms of the cluster coordinates and lattice vectors, which we solve in terms of $\<u,v\>$ and optimize over all packings in the family. The maximum occurs at $\<u,v\>\approx \<-.034789016702,+.089604971413\>$, with packing density $D\approx .778615700855$.

\bigskip 2.1\, The cluster ${\bf B}_{u,v}^9$ [See Fig. 1]

\bigskip We construct the cluster ${\bf B}_{}^1$, consisting of 1 tetrahedron. It is the convex hull of the `upper edge' $E[\<+1,+1,+1\>,\<-1,-1,+1\>]$ and `lower edge' $E[\<+1,-1,-1\>,\<-1,+1,-1\>]$.

\par\tab We attach 4 additional tetrahedra face-to-face along the `upper edge', and we call this cluster ${\bf E}_u^{4+}$. There are now 5 tetrahedra along the `upper edge', and they form the `upper half-cluster' ${\bf E}_u^{5+} = {\bf B}_{}^1\cup{\bf E}_u^{4+}$.

\par\tab We attach 4 additional tetrahedra face-to-face along the `lower edge', and we call this cluster ${\bf E}_v^{4-}$. There are now 5 tetrahedra along the `lower edge', and they form the `lower half-cluster' ${\bf E}_v^{5-} = {\bf B}_{}^1\cup{\bf E}_v^{4-}$.

\par\tab There are 9 total tetrahedra, and they form the `cluster' ${\bf B}_{u,v}^9 = {\bf E}_u^{5+}\cup{\bf E}_v^{5-}$. Clusters have 2 orientations: `positive' clusters $+{\bf B}_{u,v}^9$ and `negative' clusters $-{\bf B}_{u,v}^9$ (point reflection). The set of clusters $\{{\bf B}_{u,v}^9 : -{1\over 9}\le u\le+{1\over 9},-{1\over 9}\le v\le+{1\over 9}\}$ forms a 2-parameter family of clusters.

\bigskip 2.2\, The swivel parameters $\<u,v\>$ [See Fig. 2]

\bigskip The cluster ${\bf E}_u^{4+}$ can rotate about the `upper edge' through small angles, without overlapping with ${\bf B}_{}^1$. The `upper swivel parameter' $-{1\over 9}\le u\le+{1\over 9}$ corresponds to the $x$-coordinate of the highest vertex $q_u^+$, and the rotation angle is given by $u = \sqrt 3\sin\theta_u^{}$.

\par\tab The cluster ${\bf E}_v^{4-}$ can rotate about the `lower edge' through small angles, without overlapping with ${\bf B}_{}^1$. The `lower swivel parameter' $-{1\over 9}\le v\le+{1\over 9}$ corresponds to the $x$-coordinate of the lowest vertex $q_v^-$, and the rotation angle is given by $v = \sqrt 3\sin\theta_v^{}$.

\bigskip 2.3\, The optimized packing ${\bf P}_{u,v}^{}$ [See Fig. 3]

\bigskip The packing is made in layers. The `central' layer consists of translates of `positive' clusters, in a 2-dimensional lattice ${\bf L}_{u,v}^2 = {\bf Z}\cdot 2a+{\bf Z}\cdot 2b$ with basis $\{2a,2b\}$. (When $\<u,v\> = \<0,0\>$, this lattice is a square lattice. For our small values of $\<u,v\>$, this lattice is almost a square lattice.)

\par\tab The layers alternate between `positive' layers containing `positive' clusters, and `negative' layers containing `negative' clusters. The layer `above' the central layer consists of translates of `negative' clusters, in the offset $c-{\bf L}_{u,v}^2$. The layer `below' the central layer consists of translates of `negative' clusters, in the offset $d-{\bf L}_{u,v}^2$.

\eject\tab The `positive' clusters live in a 3-dimensional lattice ${\bf L}_{u,v}^3 = {\bf L}_{u,v}^2+{\bf Z}(c-d)$ with basis $\{2a,2b,c-d\}$, and the `negative' clusters live in the offset $c-{\bf L}_{u,v}^3 = d-{\bf L}_{u,v}^3$.
The packing consists of 2 cosets of the lattice ${\bf L}_{u,v}^3$, which correspond to the 2 orientations. 
For each cluster ${\bf B}_{u,v}^9$, we find the optimized packing and call it ${\bf P}_{u,v}^{}$.

\par\tab The packing is periodic and crystallographic, and has a large isometry group which acts transitively on clusters. There is a direct isometry between any two clusters of the same orientation, and an indirect isometry between any two clusters of opposite orientations. The direct isometry group is the translation group ${\bf L}_{u,v}^3 = {\bf Z}\cdot 2a+{\bf Z}\cdot 2b+{\bf Z}(c-d)$. The point group is (generically) trivial.

\par\tab The set of optimized packings $\{{\bf P}_{u,v}^{} : -{1\over 9}\le u\le+{1\over 9},-{1\over 9}\le v\le+{1\over 9}\}$ forms a 2-parameter family of packings.

\bigskip 2.4\, The neighbors \& separating planes of ${\bf B}_{u,v}^9$ [See Fig. 4]

\bigskip Since ${\bf B}_{u,v}^9$ is non-convex,
we consider the half-clusters ${\bf E}_u^{5+},{\bf E}_v^{5-}$ separately. Since each half-cluster is disjoint from its neighbors, the union is also disjoint from its neighbors. Therefore, we can verify that we have a packing.

\par\tab Separating planes separate the upper half-cluster ${\bf E}_u^{5+}$ from each of its neighbors: 
the plane ${\bf R}\cdot 2a+{\bf R}\cdot 2b+q_u^+$ for the upper above half-layer, 4 planes for the 6 neighbors in the lower above half-layer, 4 planes for the 8 neighbors in the upper central half-layer, 6 planes for the 8 neighbors in the lower central half-layer, and the plane ${\bf R}\cdot 2a+{\bf R}\cdot 2b+d+q_u^+$ for the upper below half-layer. ($q_u^+$ is the highest point of the upper half-cluster.)

\par\tab Separating planes separate the lower half-cluster ${\bf E}_v^{5-}$ from each of its neighbors: 
the plane ${\bf R}\cdot 2a+{\bf R}\cdot 2b+c+q_v^-$ for the lower above half-layer, 6 planes for the 8 neighbors in the upper central half-layer, 4 planes for the 8 neighbors in the lower central half-layer, 4 planes for the 6 neighbors in the upper below half-layer, and the plane ${\bf R}\cdot 2a+{\bf R}\cdot 2b+q_v^-$ for the lower below half-layer. ($q_v^-$ is the lowest point of the lower half-cluster.)

\bigskip 2.5\, The intersections \& separating planes of ${\bf B}_{u,v}^9$ [See Fig. 5]

\bigskip In the same layer, clusters touch transversely edge-to-edge, intersecting at only 1 point. Each upper half-cluster has 4 such contacts with neighboring upper half-clusters, and each lower half-cluster has 4 such contacts with neighboring lower half-clusters. These intersections give 4 equations. Each upper half-cluster has $\le 2$ such contacts with neighboring lower half-clusters, and each lower half-cluster has $\le 2$ such contacts with neighboring upper half-clusters. These intersections give $\le 2$ equations.

\par\tab  For non-intersecting skew edges $E[a,b]\ne E[c,d]+w$, there are many choices for the separating plane. We can choose the plane which contains $E[a,b]$ and parallel to $E[c,d]$ (ie: the plane through the points $a,b,b+c-d$), the plane which contains $E[c,d]$ and parallel to $E[a,b]$ (ie: the plane through the points $c,d,d+a-b$), or any parallel plane in between.

\par\tab For the transverse edge-to-edge intersection $E[a,b] = E[c,d]+w$,
the separating plane contains the union of the edges $E[a,b]\cup E[c,d]$. We solve the virtual intersection equation $sa+(1-s)b = tc+(1-t)d+w$, between lines with parameters $s,t$. For our (small) values of $\<u,v\>$, our edges intersect in their interiors (ie: $0\le s\le 1$ and $0\le t\le 1$).

\eject\tab In adjacent layers, clusters touch partially face-to-face. Each upper half-cluster has 3 such contacts with neighboring half-clusters in the above layer, and each lower half-cluster has 3 such contacts with neighboring half-clusters in the below layer. These intersections give 6 equations.

\par\tab For non-intersecting parallel faces $F[a,b,c]\ne F[d,e,f]+w$, there are many choices for the separating plane. We can choose the plane containing $F[a,b,c]$, the plane containing $F[d,e,f]$, or any parallel plane in between. 

\par\tab For the partial face-to-face intersection $F[a,b,c] = F[d,e,f]+w$, the separating plane is the plane which contains the union of the faces $F[a,b,c]\cup F[d,e,f]$. There are two ways to write the intersection equation.

\par\tab We can reduce it to a transverse edge-to-edge intersection of any non-parallel pair of edges from $\{E[b,c],E[c,a],E[a,b]\}$ and $\{E[e,f],E[f,d],E[d,e]\}$). If we choose carefully, we can find two edges which intersect in their interiors.

\par\tab Alternatively, we can reduce it to an incidence between a triangle $F[a,b,c]$ and a vertex from $\{d,e,f\}$, or vice versa. Then we solve the virtual intersection equation $sa+tb+(1-s-t)c = d+w$, with plane parameters $s,t$. If we choose carefully, we can find a point which lies inside the triangle (ie: $0\le s\le 1$ and $0\le t\le 1$).

\bigskip
\bigskip{\bf 3\, Pictures \& Equations}

\bigskip The figures are really super-figures which contain multiple sub-figures. Figs. 1-5 are floor-plan diagrams which show how the various layers fit together (ie: the top layer is placed at the top, and the bottom layer is placed at the bottom). For the purpose of perspective, we draw the coordinate axes as solid objects. For the purpose of comparison, we provide both top views (perspective vector $\<0,0,1\>$) \& side views (perspective vector $\<0,-\cos{\pi\over 10},\sin{\pi\over 10}\>$).

\medskip $\bullet$ Fig. 0 shows the orientation of the $\<x,y,z\>$ coordinate axes (for Figs. 1-5).
\par $\bullet$ Fig. 1 (2 pages) shows the cluster ${\bf B}_{u,v}^9$.
\par $\bullet$ Fig. 2 (2 pages) shows the swivel parameters $\<u,v\>$.
\par $\bullet$ Fig. 3 (2 pages) shows the packing ${\bf P}_{u,v}^{}$.
\par $\bullet$ Fig. 4 (2 pages) shows the neighbors \& separating planes (general overview).
\par $\bullet$ Fig. 5 (4 pages) shows the intersections \& separating planes (fine details).

\bigskip
\bigskip\includegraphics[height=1.5in,width=4.5in]{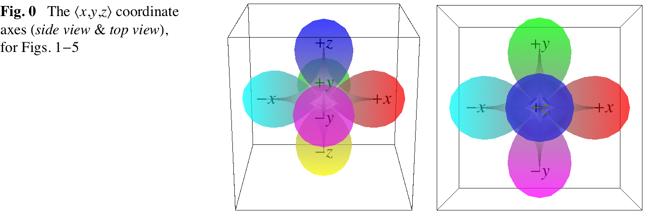}

\eject\includegraphics[height=7.666667in,width=4.5in]{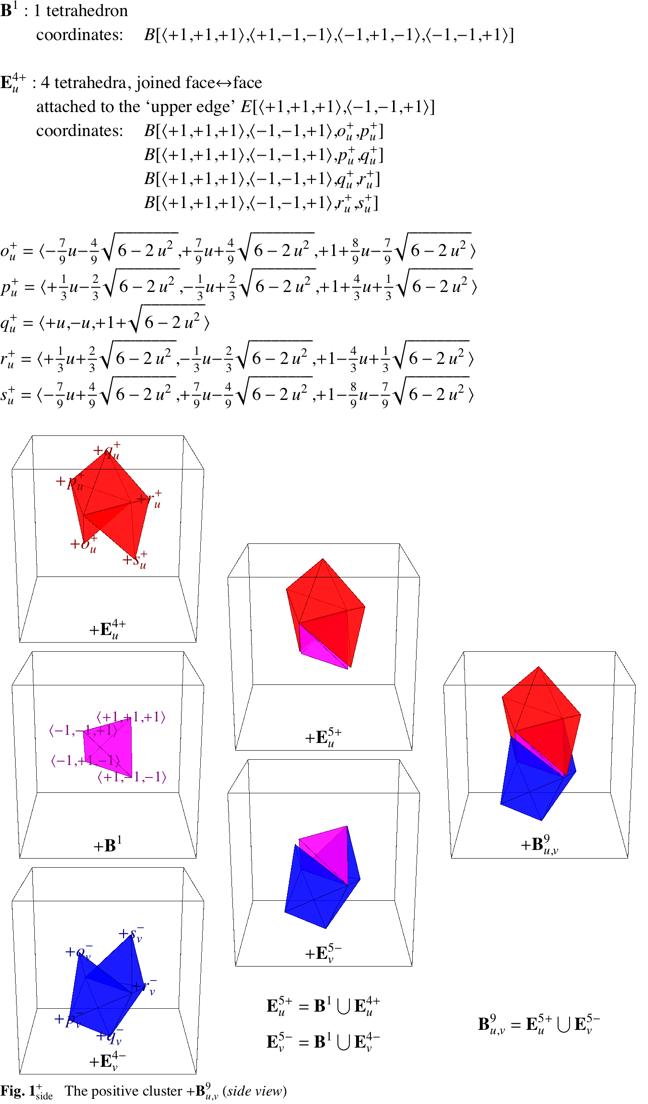}
\eject\includegraphics[height=7.666667in,width=4.5in]{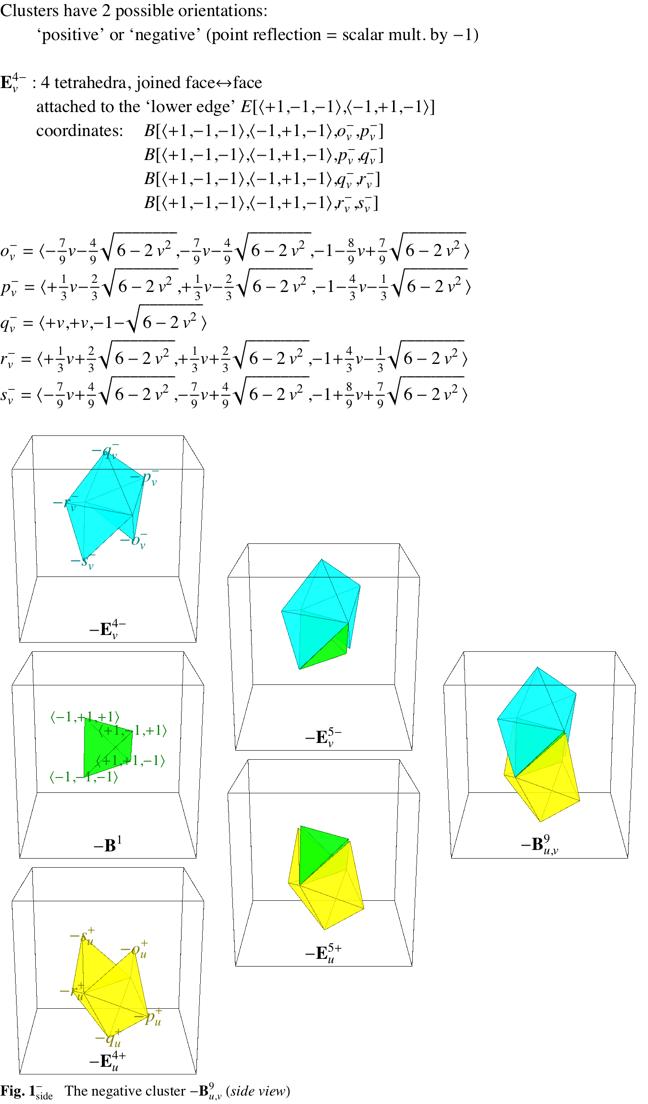}
\eject\includegraphics[height=7.666667in,width=4.5in]{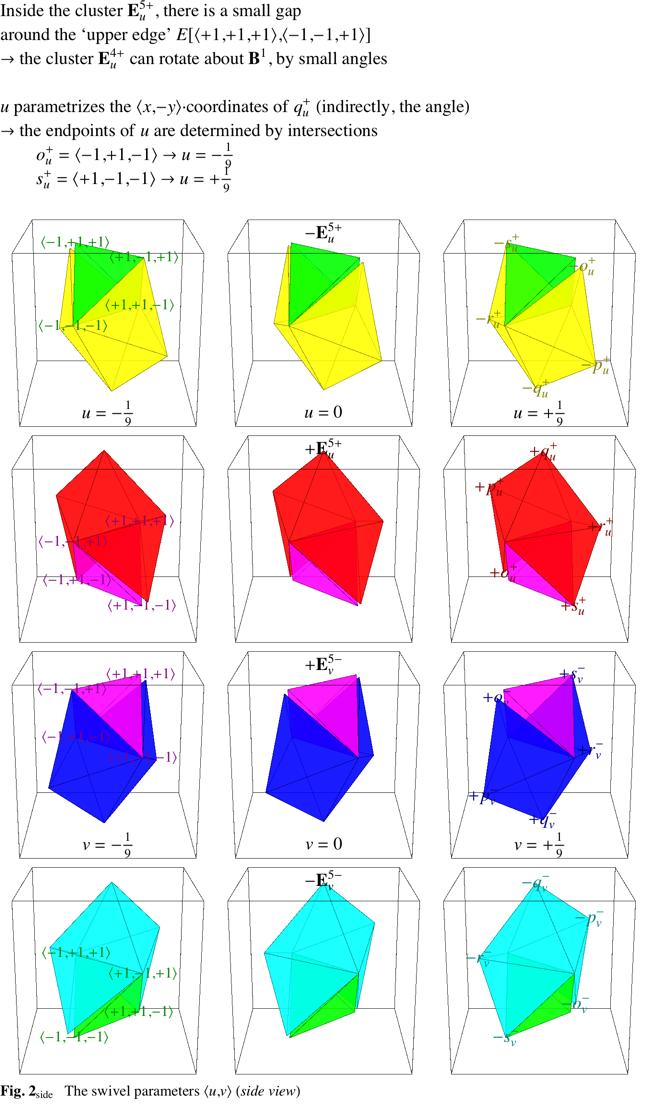}
\eject\includegraphics[height=7.666667in,width=4.5in]{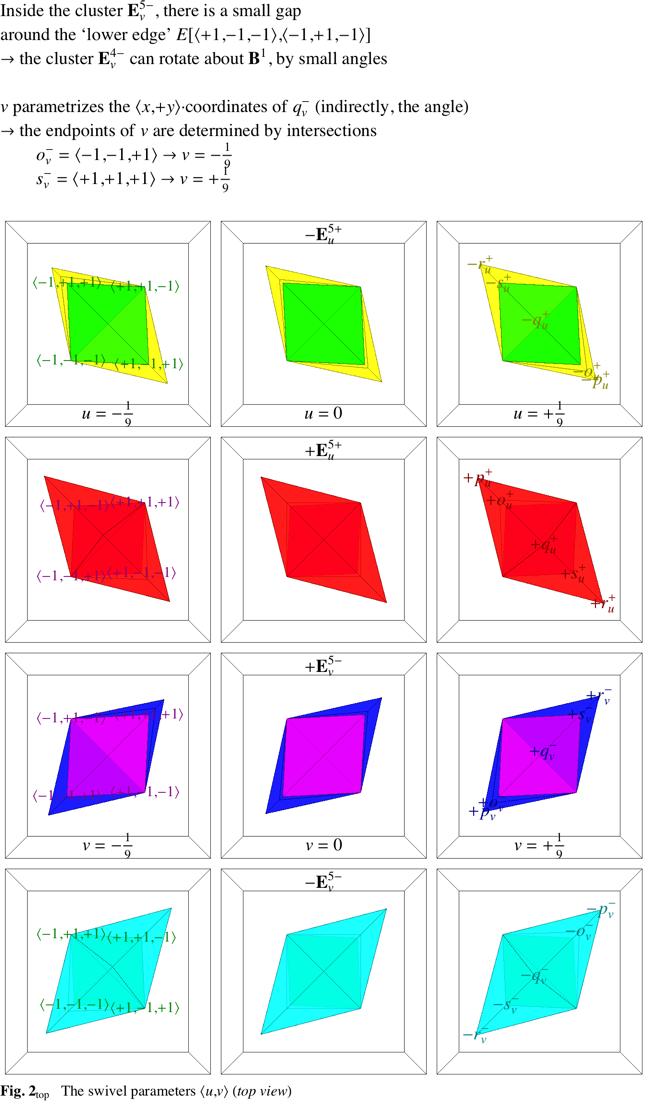}
\eject\includegraphics[height=7.666667in,width=4.5in]{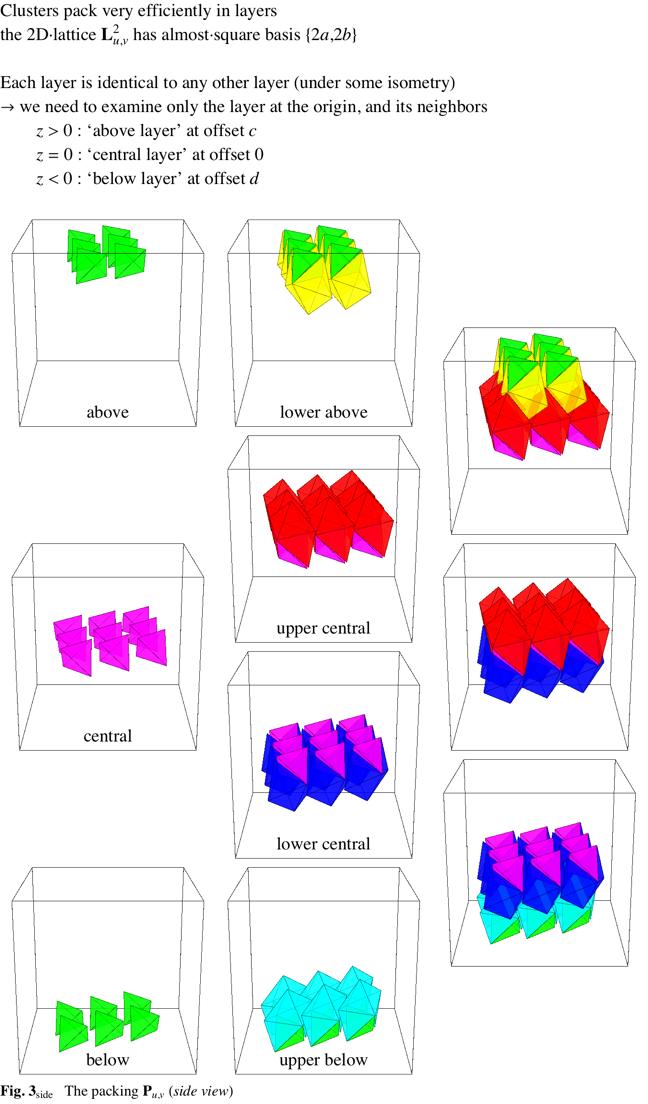}
\eject\includegraphics[height=7.666667in,width=4.5in]{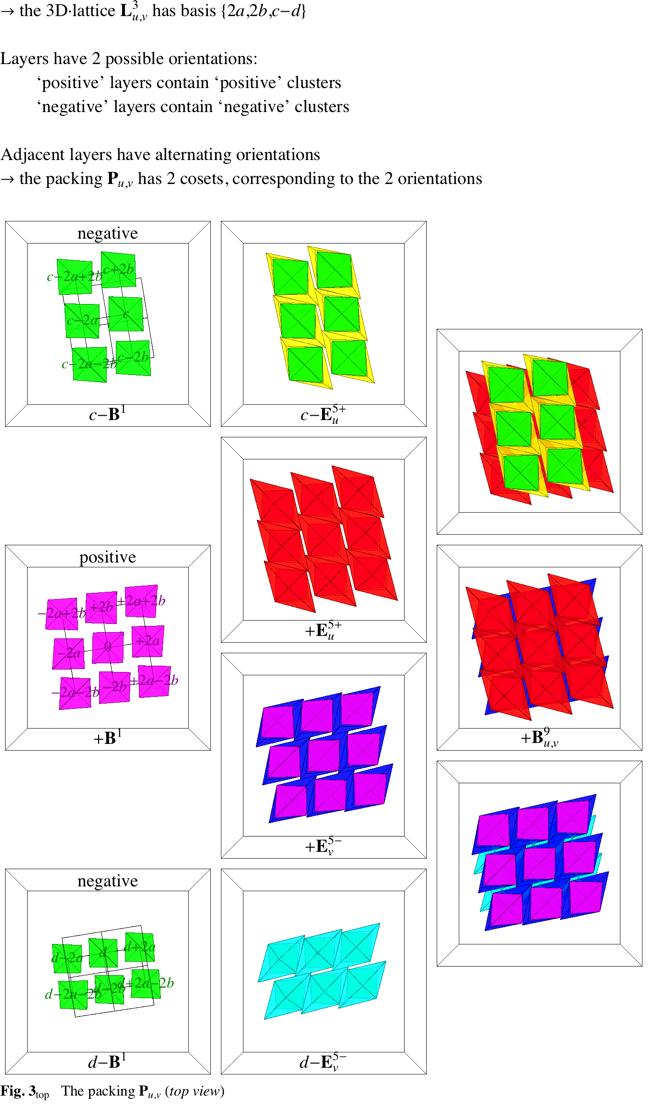}
\eject\includegraphics[height=7.666667in,width=4.5in]{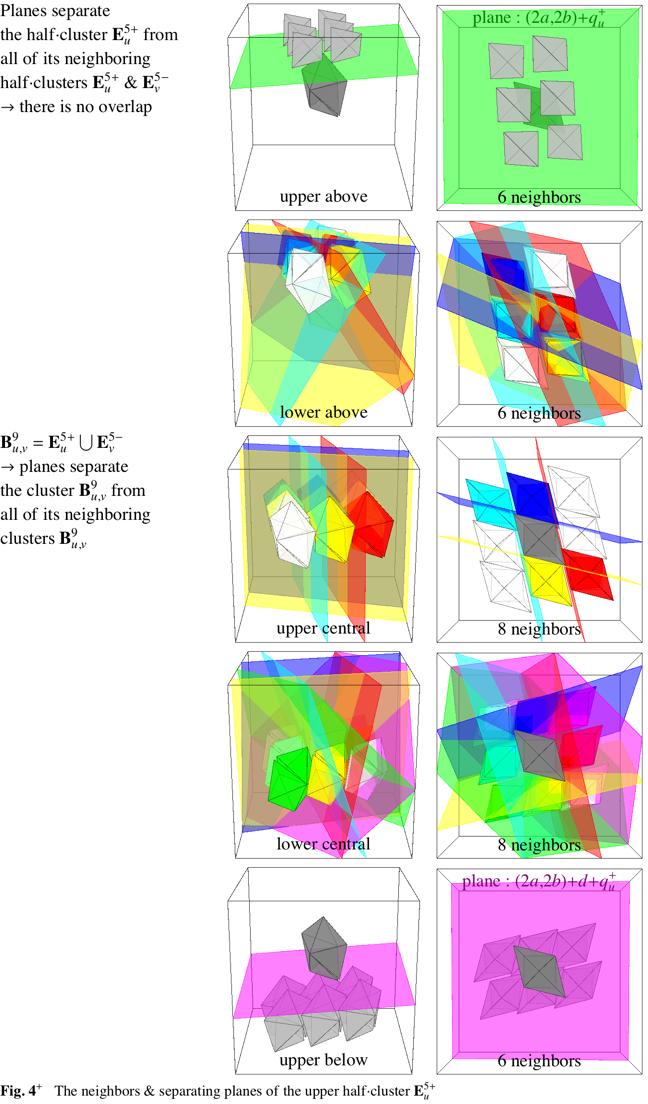}
\eject\includegraphics[height=7.666667in,width=4.5in]{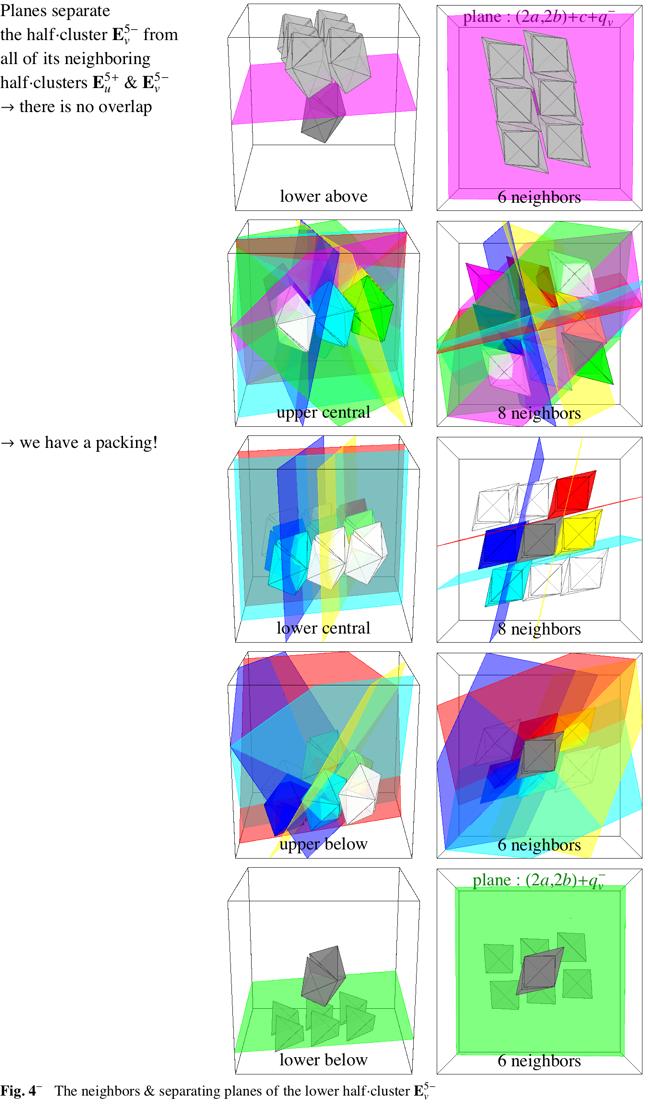}
\eject\includegraphics[height=7.666667in,width=4.5in]{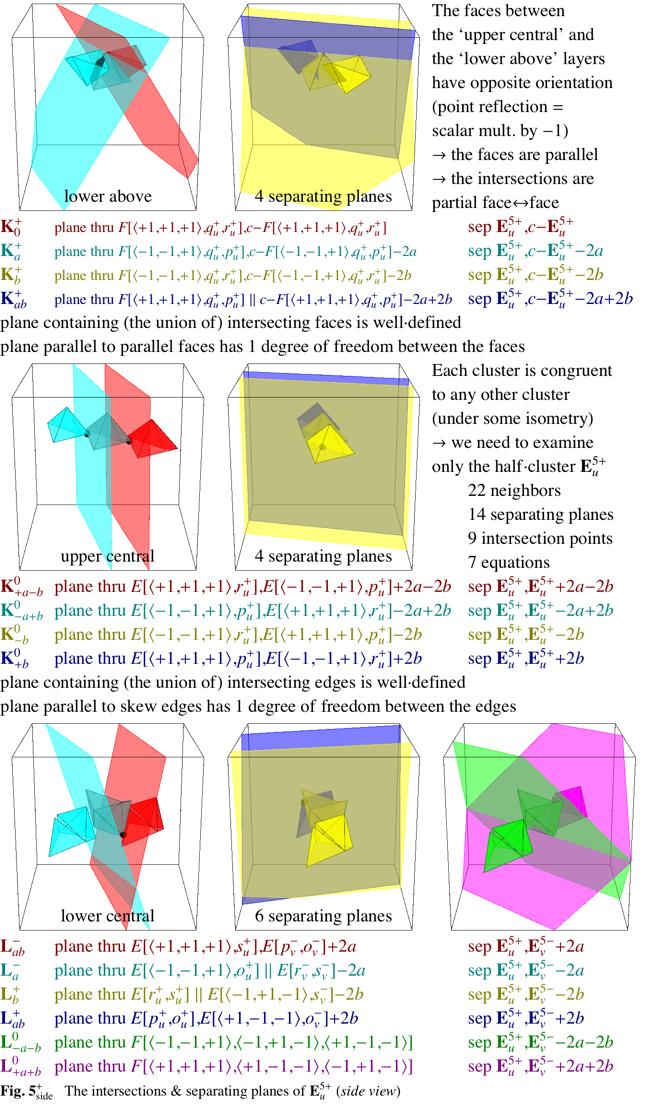}
\eject\includegraphics[height=7.666667in,width=4.5in]{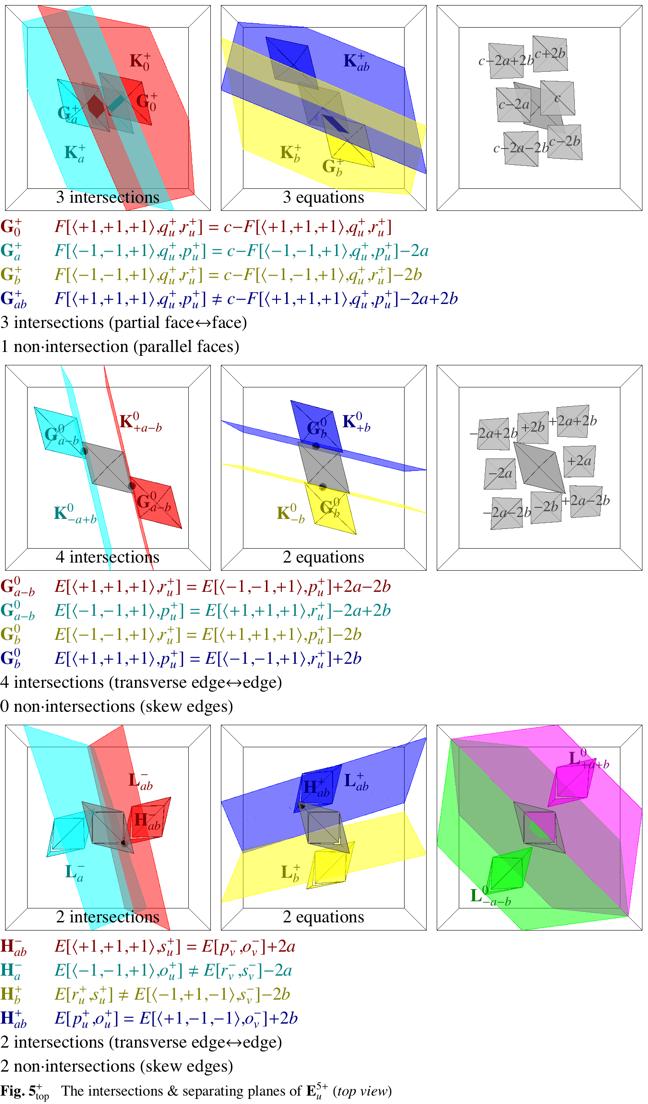}
\eject\includegraphics[height=7.666667in,width=4.5in]{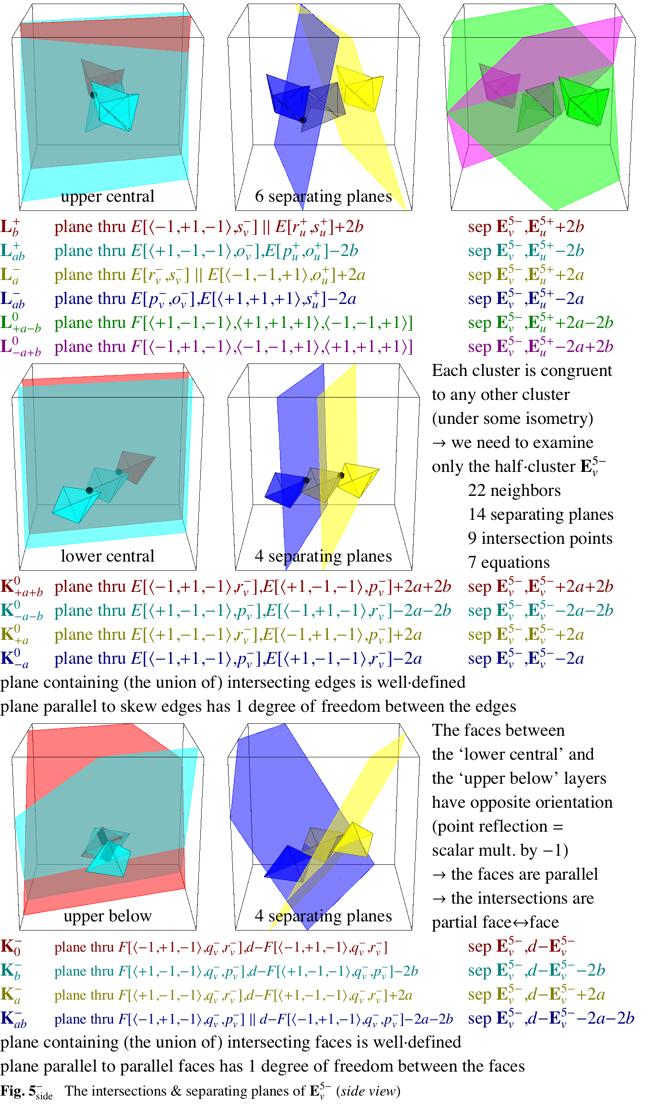}
\eject\includegraphics[height=7.666667in,width=4.5in]{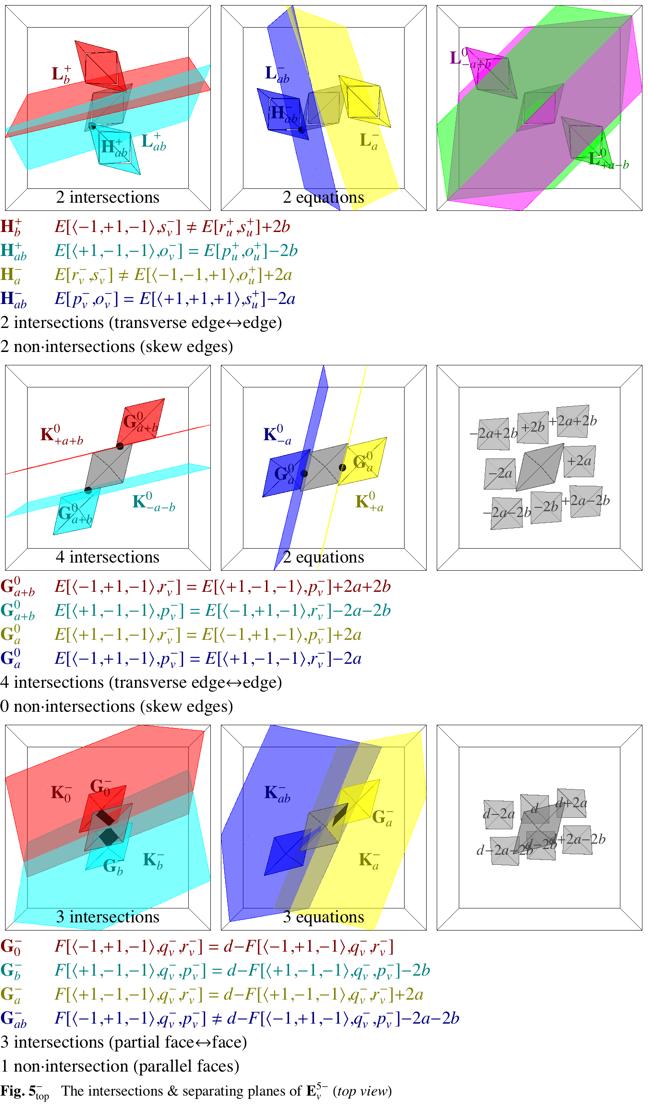}

\eject{\bf 4\, Proof by Computation}

\bigskip 4.1\, Numerical results

\bigskip We use the following procedure to compute the maximal optimized packing:

\bigskip 1. Write abstract expressions for the virtual cluster vertices
$\{o_u^+,p_u^+,q_u^+,r_u^+,s_u^+\}$,
\par $\{o_v^-,p_v^-,q_v^-,r_v^-,s_v^-\}$, in terms of the swivel parameters $\<u,v\>$. [see Fig. 1]

\bigskip 2. Write abstract equations for the virtual intersections between the cluster and its neighbors, in terms of the virtual cluster vertices (ultimately in terms of the swivel parameters), the lattice vectors $\{a,b,c,d\}$ and extra intersection parameters (edge-to-edge or point-in-plane). [see Fig. 5]

\bigskip\tabb ${\bf G} = \{{\bf G}_0^+,{\bf G}_a^+,{\bf G}_b^+,{\bf G}_0^-,{\bf G}_a^-,{\bf G}_b^-, {\bf G}_a^0,{\bf G}_b^0,{\bf G}_{a+b}^0,{\bf G}_{a-b}^0\}$.
\par\tabb ${\bf H} = \{{\bf H}_a^{},{\bf H}_b^{}\}$

\bigskip 3. Solve for the abstract lattice vectors and extra intersection parameters, in terms of the swivel parameters. Clear the common denominator.

\bigskip 4. Write abstract expressions for the lattice volume $V$ \& packing density $D$, in terms of the lattice vectors (ultimately in terms of the swivel parameters). (The fundamental domain contains 2 clusters = 18 tetrahedra. Each tetrahedron has edge length $\sqrt 8$ and volume ${8\over 3}$.)

\bigskip\tabb $V = \det[2a,2b,c-d]$
\par\tabb $D = 48/V$

\bigskip 5. Find the relative minimum of the lattice volume, with respect to the swivel parameters. (We could find only a numerical approximation, but not the minimal polynomials, because the degree is too big for our computer program.)

\bigskip\tabb $\<u,v\>\approx\<-0.034789016702,+0.089604971413\>$

\bigskip $6_{}^+$. Evaluate the vertices of the upper half-cluster.

\bigskip\tabb $o_u^+\approx\<-1.061384364770,+1.061384364770,-0.935697925928\>$
\par\tabb $p_u^+\approx\<-1.644260072209,+1.644260072209,+1.769946511051\>$
\par\tabb $q_u^+\approx\<-0.034789016702,+0.034789016702,+3.448995599962\>$
\par\tabb $r_u^+\approx\<+1.621067394407,-1.621067394407,+1.862717222257\>$
\par\tabb $s_u^+\approx\<+1.115500612974,-1.115500612974,-0.873850785124\>$

\bigskip $6_{}^-$. Evaluate the vertices of the lower half-cluster.

\bigskip\tabb $o_v^-\approx\<-1.156897066820,-1.156897066820,+0.822958681256\>$
\par\tabb $p_v^-\approx\<-1.600938143111,-1.600938143111,-1.934876528675\>$
\par\tabb $q_v^-\approx\<+0.089604971413,+0.089604971413,-3.446209700372\>$
\par\tabb $r_v^-\approx\<+1.660674790719,+1.660674790719,-1.695929938240\>$
\par\tabb $s_v^-\approx\<+1.017511555733,+1.017511555733,+0.982256408212\>$

\eject 7. Evaluate the extra intersection parameters and intersection points, in order to verify that the intersections are valid (i.e.: parameter values between 0 and 1).

\bigskip 8. Evaluate the lattice vectors.

\bigskip\tabb $a\approx\<+1.286101228632,+0.200477756509,+0.117304750804\>$
\par\tabb $b\approx\<-0.216654244567,+1.293299854677,+0.064005626691\>$
\par\tabb $c\approx\<+1.049828831839,+0.572626599359,+4.643104387948\>$
\par\tabb $d\approx\<-0.331247602500,+1.312711170924,-4.452573714586\>$

\bigskip 9. Evaluate the lattice volume \& packing density.

\bigskip\tabb $V\approx 61.647870634123$
\par\tabb $D\approx 0.778615700855$

\bigskip 4.2\, The actual density function ${\bf D}$ [see Fig. 7]

\bigskip For each cluster ${\bf B}_{u,v}^9$, we compute the optimized packing ${\bf P}_{u,v}^{}$ and packing density ${\bf D}_{u,v}^{}$. The set of all packing densities $\{{\bf D}_{u,v}^{} : -{1\over 9}\le u\le+{1\over 9},-{1\over 9}\le v\le+{1\over 9}\}$ forms the `actual' density function ${\bf D}$. The range of the graphs is $+.7767\le w\le +.7787$ for Fig. 7 and $+.7675\le w\le +.78$ for Fig. 8.

\par\tab All packings in the family satisfy the set of equations ${\bf G}$, which gives all packings the same geometric structure and lattice type. (For small values of the swivel parameters $\<u,v\>$, the intersections from the equations ${\bf G}$ remain valid.)

\par\tab However, when we try to find the maximal optimized packing using only the equations ${\bf G}$, there are 2 places where the clusters overlap. In this case, we need to impose the 2 additional equations ${\bf H}$ in order to prevent overlap.

\par\tab For a general cluster ${\bf B}_{u,v}^9$, only a subset of the equations ${\bf H}$ is necessary. So there are 4 possible optimized packings using the equations
${\bf G}$, ${\bf G}\cup\{{\bf H}_a^{}\}$, ${\bf G}\cup\{{\bf H}_b^{}\}$, ${\bf G}\cup\{{\bf H}_a^{},{\bf H}_b^{}\}$,
which gives us 4 possible packing densities ${\bf D}_{u,v}^0$, ${\bf D}_{u,v}^a$, ${\bf D}_{u,v}^b$, ${\bf D}_{u,v}^{ab}$.

\par\tab A priori, we don't know which set of equations to use for a given cluster, or which clusters satisfy a given set of equations. So we na{\"\i}vely assume that all 4 possibilities are valid for all clusters, which gives us 4 `virtual' density functions ${\bf D}_{}^0$, ${\bf D}_{}^a$, ${\bf D}_{}^b$, ${\bf D}_{}^{ab}$.

\bigskip 4.3\, The virtual density functions ${\bf D}_{}^0,{\bf D}_{}^a,{\bf D}_{}^b,{\bf D}_{}^{ab}$ [see Fig. 8]

\bigskip For each swivel parameter value $\<u,v\>$, we find the `valid' virtual density function (the set of equations which gives the densest non-overlapping packing for the given cluster). Conversely, for each virtual density function, we find the `valid' region of swivel parameter values (the set of clusters having densest non-overlapping packings which satsify the given equations).

\par\tab The 4 virtual density functions ${\bf D}_{}^0,{\bf D}_{}^a,{\bf D}_{}^b,{\bf D}_{}^{ab}$ are smooth, and their surfaces are pairwise tangent along their intersection curves ${\bf D}_{}^0\cap{\bf D}_{}^a,{\bf D}_{}^0\cap{\bf D}_{}^b,{\bf D}_{}^a\cap{\bf D}_{}^{ab},{\bf D}_{}^b\cap{\bf D}_{}^{ab}$. All 4 surfaces are mutually tangent at their common intersection point ${\bf D}_{}^0\cap{\bf D}_{}^a\cap{\bf D}_{}^b\cap{\bf D}_{}^{ab}$, which is also the common intersection point of the 4 intersection curves.

\eject\includegraphics[height=1.5in,width=4.5in]{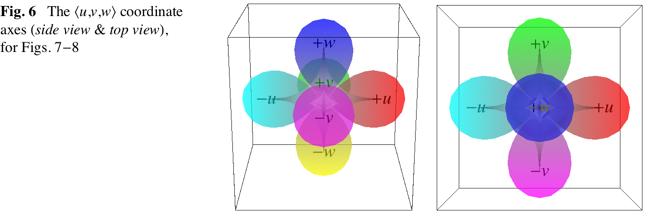}

\bigskip
\bigskip\includegraphics[height=3in,width=4.5in]{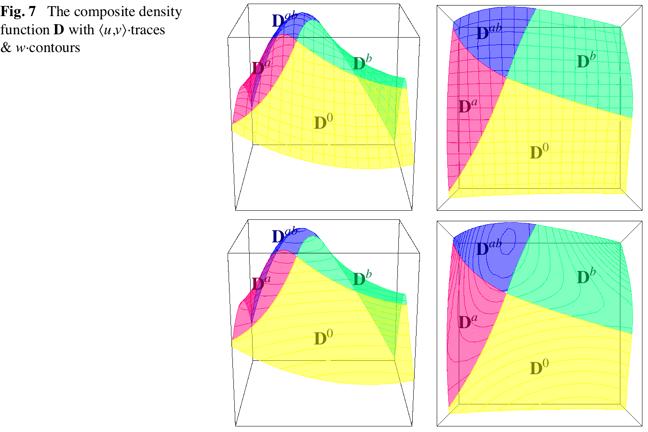}

\bigskip
\bigskip\tab The actual density function ${\bf D}$ consists of the union of the 4 valid regions of the virtual density functions, which are bounded by the intersection curves. The 4 valid regions are pairwise tangent at the intersection curves, and mutually tangent at their common intersection point. Thus, the actual density function is also smooth.

\par\tab The common intersection point occurs at

\bigskip\tabb $\<u,v,w\>\approx\<-0.037320921073,+0.033926596665,+0.778365087767\>$.

\bigskip\tab The relative maximum of the actual density function ${\bf D}$ occurs at

\bigskip\tabb $\<u,v,w\>\approx\<-0.034789016702,+0.089604971413,+0.778615700855\>$.

\bigskip $\bullet$ Figure 6 shows the orientation of the $\<u,v,w\>$ coordinate axes (for Figs. 7-8).
\par $\bullet$ Figure 7 shows the actual density function ${\bf D}$.
\par $\bullet$ Figure 8 (2 pages) shows the virtual density functions ${\bf D}_{}^0,{\bf D}_{}^a,{\bf D}_{}^b,{\bf D}_{}^{ab}$.

\eject 4.4\, The highly symmetric cluster ${\bf B}_{0,0}^9$ and packing ${\bf P}_{\rm sym}^{}$

\bigskip The highly symmetric cluster ${\bf B}_{0,0}^9$ has symmetry group ${\bf I}_{x,y,z}^{}$.

\bigskip\tabb ${\bf I}_{x,y,z}^{} = \{I,P,P_{}^2,P_{}^3,Q,QP,QP_{}^2,QP_{}^3\}$,
$P = \scriptsize\left[\matrix{0&-1&0\cr +1&0&0\cr 0&0&-1}\right]$,
$Q = \scriptsize\left[\matrix{+1&0&0\cr 0&-1&0\cr 0&0&-1}\right]$

\bigskip\tab The highly symmetric packing ${\bf P}_{\rm sym}$ also has symmetry group ${\bf I}_{x,y,z}^{}$. The lattice has a square basis, with lattice vectors of the form

\bigskip\tabb $a = \<+2i,+2j,0\>$
\par\tabb $b = \<-2j,+2i,0\>$
\par\tabb $c = \<+i,+j,+k\>$
\par\tabb $d = \<-j,+i,-k\>$,

\bigskip\tabb $i = {1\over 71}(-168+106\sqrt 6)\approx 1.290787503310$
\par\tabb $j = {1\over 71}(-88+42\sqrt 6)\approx 0.209557312632$
\par\tabb $k = {1\over 71}(-262+238\sqrt 6)\approx 4.520824771583$

\bigskip and lattice volume \& packing density

\bigskip\tabb $V = {1\over 357911}(-730200320+307139840\sqrt 6)\approx 61.846569901642$
\par\tabb $D = {1\over 4477040}(1711407+719859\sqrt 6)\approx .776114181859$.

\bigskip The highly symmetric packing ${\bf P}_{\rm sym}$ is {\it not} a member of the family $\{{\bf P}_{u,v}^{}\}$, because adjacent layers have only 2 face-to-face intersections.

\par\tab The optimized packing ${\bf P}_{0,0}^{}$ has trivial symmetry group. As with all members of the family, adjacent layers have 3 face-to-face intersections. In order to get this additional intersection, we need to shift adjacent layers by a small distance (relative to the lattice), which destroys symmetry.

\bigskip 4.5 Isometries of the cluster ${\bf B}_{u,v}^9$ and packing ${\bf P}_{u,v}^{}$ [see Fig. 9]

\bigskip We can apply the isometries ${\bf I}_{x,y,z}$ in several ways:

\bigskip 1. Apply the isometries to the optimized packing ${\bf P}_{u,v}^{}$ (which also applies the isometries to the clusters and intersection equations).
The 8 clusters ${\bf I}_{x,y,z}\cdot{\bf B}_{u,v}^9$ are congruent to each other, the 8 optimized packings ${\bf I}_{x,y,z}\cdot{\bf P}_{u,v}^{}$ are congruent to each other, and the 8 optimized packing densities are equal to each other.

\par\tab Each optimized packing corresponds to a point of an actual density function. The 8 different (but analogous) points have the same $w$-coordinate (packing density), but each point lives in a different (but congruent) surface (actual density functions ${\bf I}_{u,v,w}\cdot{\bf D}$). (The isometries ${\bf I}_{u,v,w}$ are similar to the isometries ${\bf I}_{x,y,z}$, except that the $w$-coordinate (packing density) is always positive.)

\eject\includegraphics[height=4.166667in,width=4in]{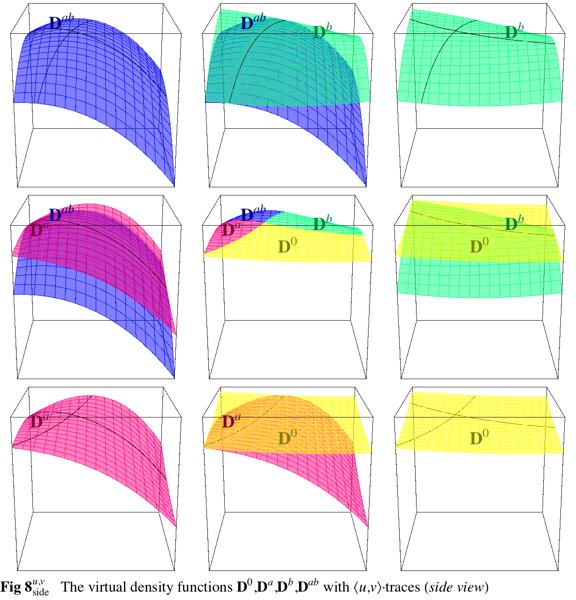}

\bigskip
\bigskip\tabb ${\bf I}_{u,v,w}^{} = \{I,R,R_{}^2,R_{}^3,S,SR,SR_{}^2,SR_{}^3\}$,
$R = \scriptsize\left[\matrix{0&-1&0\cr +1&0&0\cr 0&0&+1}\right]$,
$S = \scriptsize\left[\matrix{+1&0&0\cr 0&-1&0\cr 0&0&+1}\right]$

\bigskip 2. Apply the isometries to the cluster ${\bf B}_{u,v}^9$, but fix the intersection equations ${\bf G}\cup{\bf H}$. The 8 clusters ${\bf I}_{x,y,z}\cdot{\bf B}_{u,v}^9$ are identical to

\bigskip\tabb $\{{\bf B}_{+u,+v}^9,{\bf B}_{-v,+u}^9,{\bf B}_{-u,-v}^9,{\bf B}_{+v,-u}^9,{\bf B}_{+v,+u}^9,{\bf B}_{+u,-v}^9,{\bf B}_{-v,-u}^9,{\bf B}_{-u,+v}^9\}$

\bigskip and congruent to each other. The 8 optimized packings

\bigskip\tabb $\{{\bf P}_{+u,+v}^{},{\bf P}_{-v,+u}^{},{\bf P}_{-u,-v}^{},{\bf P}_{+v,-u}^{},{\bf P}_{+v,+u}^{},{\bf P}_{+u,-v}^{},{\bf P}_{-v,-u}^{},{\bf P}_{-u,+v}^{}\}$

\bigskip are {\it not} congruent to each other, and the 8 optimized packing densities

\bigskip\tabb $\{{\bf D}_{+u,+v}^{},{\bf D}_{-v,+u}^{},{\bf D}_{-u,-v}^{},{\bf D}_{+v,-u}^{},{\bf D}_{+v,+u}^{},{\bf D}_{+u,-v}^{},{\bf D}_{-v,-u}^{},{\bf D}_{-u,+v}^{}\}$

\bigskip are {\it not} equal to each other.

\par\tab Each optimized packing corresponds to a point of an actual density function. All 8 points live in the same surface (actual density function ${\bf D}$), but each point has a different $w$-coordinate (packing density).

\eject\includegraphics[height=4.166667in,width=4in]{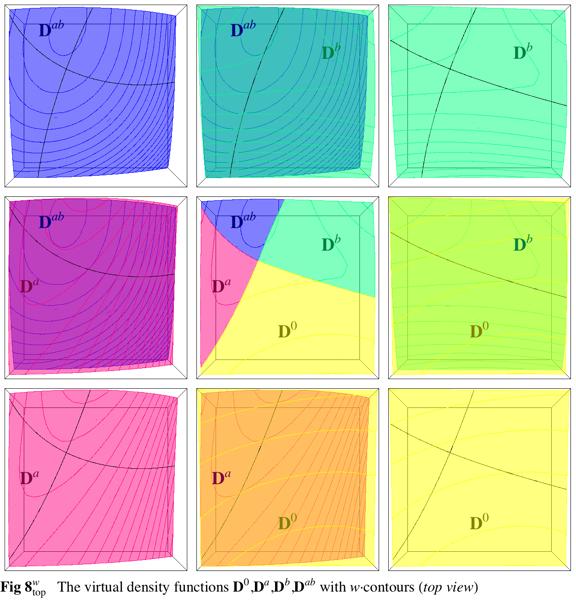}

\bigskip
\bigskip 3. Apply the isometries to the intersection equations ${\bf G}\cup{\bf H}$, but fix the cluster ${\bf B}_{u,v}^9$. This is equivalent to fixing the intersection equations and applying the inverse isometries to the cluster. Since ${\bf I}_{x,y,z}^{}$ is a group, the set of packings \& densities is identical.

\par\tab You could also apply the isometries individually to each vector quantity in the equations. For the virtual cluster vertices $\{o_u^+,p_u^+,q_u^+,r_u^+,s_u^+\}$, $\{o_v^-,p_v^-,q_v^-,r_v^-,s_v^-\}$, use the highly symmetric cluster ${\bf B}_{0,0}^{}$ as a paradigm. For the abstract lattice vectors $\{a,b,c,d\}$, use the highly symmetric packing ${\bf P}_{\rm sym}^{}$ as a paradigm.

\par\tab In any case, we don't get any additional packings by considering additional isometries of the cluster, intersection equations or packings. Thus, it suffices to optimize over a single generic cluster ${\bf B}_{u,v}^9$, a single set of equations ${\bf G}\cup{\bf H}$, and a single actual density function ${\bf D}$.

\par\tab The 8 congruences of the cluster ${\bf B}_{u,v}^9$ are symmetric (about $\langle u,v\rangle = \langle 0,0\rangle$) as a whole, so the 8 sets of intersection equations ${\bf G}\cup{\bf H}$ are symmetric as a whole, and the 8 congruences of the density function ${\bf D}$ are symmetric as a whole. However, any single set of intersection equations is individually asymmetric, so any single density function is individually asymmetric.

\bigskip $\bullet$ Figure 9 (2 pages) shows the isometries of the cluster ${\bf B}_{u,v}^9$ and packing ${\bf P}_{u,v}^{}$.

\eject\includegraphics[height=7.333333in,width=4.666667in]{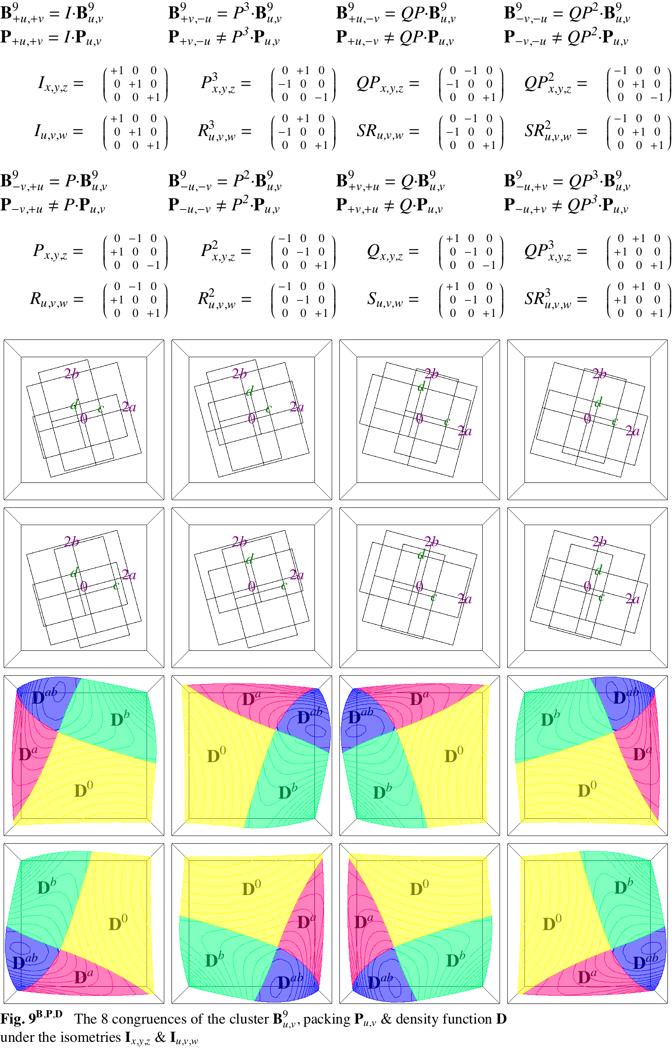}
\eject\includegraphics[height=7.166667in,width=4.666667in]{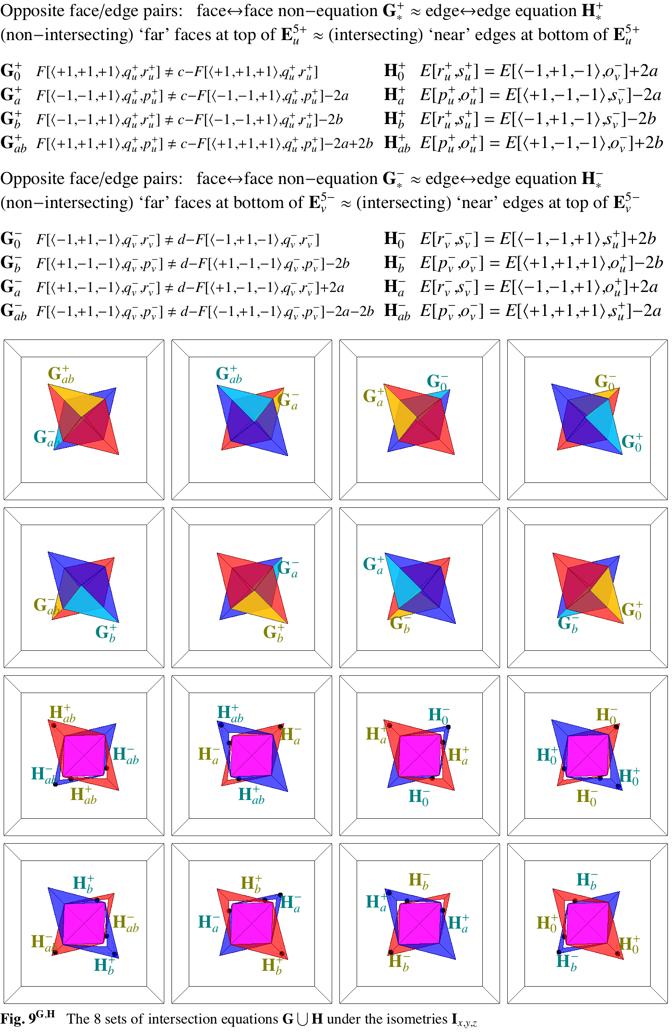}

\eject 4.6\, Remarks: Dense Packing

\bigskip This construction gives a dense packing, for several reasons:

\medskip 1. The half-cluster ${\bf E}_{}^5$ (5 tetrahedra joined face-to-face along a common edge) is very dense locally. The local density (fraction of solid angle around the edge) is $D\approx .979566380077$. The relative density of the half-cluster inside its convex hull is $D = {405\over 416}\approx .973557692308$.

\par 2. Opposite edges of a regular tetrahedron are skew perpendicular.
Inside the cluster ${\bf B}_{u,v}^9$, the two half-clusters are oriented perpendicularly. So the clusters fit together very nicely in layers with an almost-square basis.

\par 3. Between adjacent layers, the half-clusters and their neighbors are point reflections (scalar mult. by $-1$). Therefore, the faces are parallel, and the intersections are partial face-to-face.

\par 4. For small values of $\<u,v\>\approx\<0,0\>$, the clusters have an extra almost-symmetry. Because the layers have almost-square basis, point reflection (scalar mult. by $-1$) is almost-equivalent to rotation by ${\pi\over 2}$ radians about the $z$-axis. Thus, the layers fit together very well to form a dense packing.

\bigskip
\bigskip{\bf 5\, Conclusion}

\bigskip For the ${\bf B}_{}^9$ cluster, if we rotate the 8 non-central tetrahedra independently about the central tetrahedron ${\bf B}_{}^1$, we can generalize to an 8-parameter family ${\bf B}_{a,b,c,d,e,f,g,h}^9 = {\bf B}_{}^1\cup{\bf E}_{a,b,c,d}^{4+}\cup{\bf E}_{e,f,g,h}^{4-}$.

\par\tab Intuitively, the maximal packing occurs when there are no gaps among the 4 upper tetrahedra, and no gaps among the 4 lower tetrahedra (ie: the special case of our 2-parameter family ${\bf B}_{u,v}^9 = {\bf B}_{u,u,u,u,v,v,v,v}^9$, where the non-central tetrahedra form ${\bf E}_u^{4+} = {\bf E}_{u,u,u,u}^{4+}$ and ${\bf E}_v^{4-} = {\bf E}_{v,v,v,v}^{4-}$).

\par\tab From a geometric perspective, it's easy to see why. If you start with a gapless cluster, and try to open a gap somewhere, you will only decrease the packing density. The distance between neighbors in the same layer will increase,
and the perpendicular distance between layers will not decrease. Therefore, generalizing will not help us.

\par\tab I found packings for other clusters of tetrahedra ${\bf B}_{17}^{}$, ${\bf V}_{20}^{}$, ${\bf E}_{15}^{}$, ${\bf B}_{13}^{}$, ${\bf B}_9^{}$ with packing density $D\ge{7\over 10}$, but the winner is ${\bf B}_9^{}$.

\par\tab However, it's possible that other arrangements of tetrahedra yet undiscovered could give a denser packing...

\bigskip
\bigskip{\footnotesize{\bf Acknowledgements}\, A huge thanks to my advisor Jeffrey Lagarias, for suggesting this problem, guiding me along the entire paper writing and submission process, and being such a meticulous advisor as well as a caring mentor! Really, without Jeff, none of this would have been possible!!! Thanks to Martin Henk for very generously sharing his fabulous computer program with me, and helping me to use it! Thanks to Pedagoguery Software Inc for generously providing many 3D models of regular tetrahedra, at cost. Many thanks to Tom Hales for writing a computer program to independently verify this packing, on very short notice, especially at a time when most people were incredulous about my result! Thanks to the reviewers for their very thorough, thoughtful and helpful comments about how to write in general, and improve this paper in particular. Many thanks to Richard Pollack for truly appreciating and supporting the artistic aspects of this project, and encouraging me to continue, even when I was feeling really disgruntled with the whole process!!! Finally, a huge thanks to everyone at Springer and VTeX for allowing me multiple revisions to get the layout exactly perfect, dealing with my quirks and eccentricities as an artist, and processing many graphics files and corrections!!

}\eject{\bf References}

\bigskip\footnotesize 1. '$\!{A}\rho\iota\sigma\tau{o}\tau\acute\epsilon\lambda\eta\varsigma$: $\Pi\epsilon\rho\grave\iota$ '$\!{O}\upsilon\rho\alpha\nu{o}\tilde\upsilon$.
Translation: Boethius: Aristotel\=es. De Caelo.
Translation: Guthrie, W.K.C.: Aristotle. On heavenly bodies. Leob Classical Library, vol. 338. Harvard University Press, vol. 6 (1986)

\par 2. Betke, U., Henk, M.: (FORTRAN computer program) lattice\_packing.f (1999)

\par 3. Betke, U., Henk, M.: Densest lattice packings of 3-polytopes. Comput. Geom. {\bf 16}(3), 157-186 (2000)

\par 4. Conway, J.H., Torquato, S.: Packing, tiling, and covering with tetrahedra. Proc. Natl. Acad. Sci. U.S.A. {\bf 103}, 10612-10617 (2006)

\par 5. Gr\"omer, H.: \"Uber die dichteste gitterf\"ormige lagerung kongruenter Tetraeder. Mon. Math. {\bf 66}, 12-15 (1962)

\par 6. Hales, T.C.: A proof of the Kepler conjecture. Ann. Math. {\bf 162}, 1065-1185 (2005)

\par 7. Hales, T.C., Ferguson, S.P.: Historical overview of the Kepler conjecture; A formulation of the Kepler conjecture; Sphere packings III: Extremal cases; Sphere packings IV: Detailed bounds; Sphere packings V: Pentahedral prisms; Sphere packings VI: Tame graphs and linear programs.
Discrete Comput. Geom. {\bf 36}(1), 5-265 (2006)

\par 8. Hilbert, D.C.: Mathematische Probleme. Nachr. Ges. Wiss. G\"ott., Math. Phys. Kl. {\bf 3}, 253-297 (1900).
Translation: Hilbert, D.C.: Mathematical problems. Bull. Am. Math. Soc. {\bf 8}, 437-479 (1902)

\par 9. Hoylman, D.J.: The densest lattice packing of tetrahedra. Bull. Am. Math. Soc. {\bf 76}, 135-137 (1970)

\par 10. Hurley, A.C.: Some helical structures generated by reflexions. Aust. J. Phys. {\bf 38}(3), 299-310 (1985)

\par 11. Minkowski, H.: Geometrie der Zahlen. Teubner, Leipzig (1896).
Reprint: Minkowski, H.: Geometrie der Zahlen. Chelsea (1953)

\par 12. Senechal, M.: Which tetrahedra fill space? Math. Mag. {\bf 54}(5), 227-243 (1981)

\par 13. Struik, D.J.: De impletione loci. Nieuw Arch. Wiskd. {\bf 15}, 121-134 (1925)

\end{document}